\documentclass[11pt]{article}
\usepackage{amsthm,amsfonts,amsmath,amssymb}
\usepackage[
    colorlinks=true,
    linkcolor=blue,
    citecolor=blue,
    urlcolor=blue
]{hyperref}
\newtheorem{thm}{Theorem}[section]
\newtheorem{prop}[thm]{Proposition}
\newtheorem{lem}[thm]{Lemma}
\newtheorem{defi}[thm]{Definition}
\newtheorem{rem}[thm]{Remark}
\newtheorem{conj}[thm]{Conjecture}
\newtheorem{cor}[thm]{Corollary}

\numberwithin{equation}{section}

\begin{document}
\title{Quadratic Elements in Finite $p$-Groups and Applications to Oliver's $p$-Group Conjecture}
\author{Jingjing Duan and Lijian An\thanks{Corresponding author. e-mail: anlj@sxnu.edu.cn }\\
Department of Mathematics,
Shanxi Normal University\\
Taiyuan, Shanxi 030031, P. R. China\\
 }

\date{}
\maketitle

\begin{abstract}
We establish a property of quadratic elements in faithful modules for
finite $p$-groups. Let $p$ be an odd prime, and suppose that $a$ is a
quadratic element of an abelian subgroup $L$ normalized by $b$.
Writing $a_0=a$ and $a_i=[a_{i-1},b]$, we prove that if, for
some $r$ with $2\leq r\leq(p+1)/2$, one has
$a_{r-1}\neq1$ and $a_r=1$, then $a_{r-1}$ is quadratic.
As an application, we verify the Weakly Closed Conjecture for finite
$p$-groups of nilpotency class at most $p+2$, and hence establish
the Quadratic Conjecture for these groups. Consequently, Oliver's $p$-group conjecture
holds for a finite $p$-group $S$ whenever
$S/\mathfrak X(S)$ has nilpotency class at most $p+2$.
\end{abstract}



\section{Introduction}

The study of quadratic elements in finite groups can be traced
back to Glauberman's work on unipotent linear groups
\cite{Glauberman1972}.

Throughout, $p$ is an odd prime unless stated otherwise, and all
modules are finite-dimensional right modules over $\mathbb F_p$.
We use the
commutator convention
\[
[x,y]=x^{-1}y^{-1}xy=x^{-1}x^y,\qquad x^y=y^{-1}xy.
\]
For subgroups $A,B\leq G$, the subgroup $[A,B]$ is generated by
all $[a,b]$ with $a\in A$ and $b\in B$. Iterated commutators are
left-normed, and $C_G(A,B)=C_G(A)\cap C_G(B)$.
For a right $\mathbb F_pG$-module $V$, put $[v,g]=v(g-1)$ and
let $[V,H]$ be the subspace spanned by all $[v,h]$, where
$v\in V$ and $h\in H$. Write $[V,H,K]=[[V,H],K]$.
The upper central series is defined by
\[
Z_0(G)=1,\qquad
Z_{i+1}(G)/Z_i(G)=Z(G/Z_i(G))\quad(i\geq0).
\]

Let $G$ be a finite $p$-group and $V$ a faithful
$\mathbb{F}_pG$-module. A nonidentity element $g\in G$ is called
quadratic on $V$ if
\[
V(g-1)^2=0.
\]
A nontrivial elementary abelian subgroup $E\leq G$ is called
\emph{quadratic on $V$} if $[V,E,E]=0$.
Quadratic elements played an important role in the
study of Oliver's $p$-group conjecture. In 2008, Green,
H\'ethelyi and Lilienthal reformulated Oliver's conjecture
in terms of $F$-modules \cite{GHL}.
A key ingredient in their work is the following property
of quadratic elements:
\begin{prop}[{\cite[Lemma 4.1]{GHL}}]\label{GHLlemma}
Let $p$ be an odd prime, $G$ a nontrivial finite
$p$-group, and $V$ a faithful $\mathbb{F}_pG$-module.
Suppose that $A,B\in G$ satisfy
\[
1\neq [A,B]\in C_G(A,B).
\]
If $[A,B]$ is nonquadratic on $V$, then both
$A$ and $B$ are nonquadratic on $V$.
\end{prop}
 Using this property, they
proved that if $G$ has nilpotency class at most two and
contains a quadratic element, then $\Omega_1(Z(G))$ also
contains a quadratic element. This result
played a key role in their verification of Oliver's
conjecture for groups $S$ such that $S/\mathfrak{X}(S)$
has nilpotency class at most two.

In subsequent work, Green, H\'ethelyi and Mazza verified Oliver's
conjecture when the relevant quotient group has nilpotency class at
most four \cite{GHM2010}.


Building on these developments, Green, H\'ethelyi and Mazza
introduced a stronger form of Oliver's conjecture in 2011, whose
module-theoretic formulation is known as the \emph{Quadratic
Conjecture} \cite{GHM2}. It asserts that if $V$ is a faithful
$F$-module for a finite $p$-group $G$, then $\Omega_1(Z(G))$
contains a nonidentity element that acts quadratically on $V$.

In this paper, we give a generalization of Proposition \ref{GHLlemma}:

\begin{thm}\label{thm:long}
Let $p$ be an odd prime, $G$ a nontrivial finite
$p$-group, and $V$ a faithful $\mathbb{F}_pG$-module.
Suppose that $L$ is an abelian subgroup of $G$.
For $a\in L$ and $b\in N_G(L)$, define a sequence of elements of $L$
recursively by
\[
a_0=a,\qquad a_{i}=[a_{i-1},b]
\quad (i\geq 1).
\]
Suppose that $a$ is quadratic on $V$ and that, for some
integer $r$ with
\[
2\leq r\leq \frac{p+1}{2},
\]
the following conditions hold:
\[
a_{r-1}\neq 1,\qquad a_r=1.
\]
Then $a_{r-1}$ is also quadratic on $V$.
\end{thm}

Theorem~\ref{thm:long} implies Proposition~\ref{GHLlemma}.
Indeed, write $C=[A,B]$ and suppose that
$1\neq C\in C_G(A,B)$. If $A$ is quadratic, take
$a=A$, $b=B$ and $L=\langle A,C\rangle$. Since
$A^B=AC$ and $C^B=C$, the element $B$ normalizes $L$;
Theorem~\ref{thm:long} with $r=2$ shows that $C$ is quadratic. Thus if $C$ is nonquadratic,
then $A$ is nonquadratic, and $B$ is nonquadratic by the same argument.

Following \cite{GHM2}, an abelian subgroup $A\leq G$ is called
\emph{weakly closed in $C_G(A)$ with respect to $G$} if
$[A,A^g]\neq1$ for every $g\in G$ such that $A^g\neq A$.
We abbreviate this to \emph{weakly closed} throughout the paper.
The Weakly Closed Conjecture of Green, H\'ethelyi and Mazza
\cite[Conjecture~4.6]{GHM2} asserts that if a faithful
$\mathbb F_pG$-module has a nontrivial weakly closed quadratic
elementary abelian subgroup, then $\Omega_1(Z(G))$ contains a
quadratic element. This assertion is false in general:
Green and Lynd constructed a group of order $3^{49}$ and a faithful
module for which it fails, although that group satisfies the
Quadratic Conjecture \cite[Theorem~1.1]{GreenLynd2013}.
Our next result verifies the assertion under a nilpotency-class bound.

\begin{thm}\label{thm:weakly-class}
Let $p$ be an odd prime, let $G$ be a finite $p$-group of nilpotency
class at most $p+2$, and let $V$ be a faithful $\mathbb F_pG$-module.
If $G$ contains a nontrivial weakly closed elementary abelian
subgroup $E$ such that $[V,E,E]=0$, then $\Omega_1(Z(G))$
contains a quadratic element.
\end{thm}

By \cite[Proposition~4.5]{GHM2}, every faithful $F$-module has
a nontrivial weakly closed quadratic offender. Thus
Theorem~\ref{thm:weakly-class} has the following consequence.

\begin{cor}\label{thm:class}
Let $p$ be an odd prime and $G$ a finite $p$-group. If the
nilpotency class of $G$ is at most $p+2$, then $G$
satisfies the Quadratic Conjecture.
\end{cor}

Using the module-theoretic criterion of Green, H\'ethelyi and
Lilienthal, we also obtain the following consequence.

\begin{cor}\label{thm:Oclass}
Let $p$ be an odd prime and $S$ a finite $p$-group. If the
nilpotency class of $S/\mathfrak{X}(S)$ is at most
$p+2$, then $S$ satisfies Oliver's $p$-group conjecture.
\end{cor}

Green, H\'ethelyi and Mazza established the Quadratic Conjecture for
finite $p$-groups of nilpotency class at most $4$
\cite{GHM2010,GHM2}. The bound $p+2$ in
Corollary~\ref{thm:class} is strictly larger. In particular, any
counterexample to the Quadratic Conjecture must have nilpotency class
at least $p+3$, giving a lower bound that grows linearly with $p$.

Regarding Oliver's $p$-group conjecture, Lynd \cite{Lynd2013}
proved that it holds whenever $p\geq5$ and the nilpotency class of
$S/\mathfrak{X}(S)$ is at most $\log_2(p-2)+1$. For every odd
prime $p\geq5$, the bound $p+2$ in Corollary~\ref{thm:Oclass}
is strictly larger. Thus Corollary~\ref{thm:Oclass} improves upon
Lynd's nilpotency-class bound.

The paper is organized as follows. In Section~2, we prove
Theorem~\ref{thm:long} and establish consequences for the upper
central series. In Section~3, we recall the relevant background
about Oliver's $p$-group conjecture and weakly closed quadratic
offenders. In Section~4, we use relative lower central series to
prove Theorem~\ref{thm:weakly-class} and
Corollaries~\ref{thm:class} and~\ref{thm:Oclass}.

\section{Proof of Theorem \ref{thm:long}}

For each $g\in G$, let $R_g$ denote the linear
transformation induced by $g$ on $V$, so that
\[
v^g=vR_g,\qquad R_{gh}=R_gR_h.
\]
We write $I$ for the identity transformation on $V$ and $I_K$
for the identity map on $K$ whenever $K$ is regarded as a vector
space. Thus a nonidentity element $g\in G$ is quadratic on $V$ if and
only if
\[
(R_g-I)^2=0.
\]
\begin{lem}\label{lem:quadratic}
Let $G$ be a finite group and let $V$ be a faithful
$\mathbb F_pG$-module. If $g\in G$ is quadratic on $V$,
then $g^p=1$. Moreover, every conjugate of $g$
is also quadratic on $V$.
\end{lem}

\begin{proof}
Put $X=R_g-I$. Since $X^2=0$, we have
$R_g^p=(I+X)^p=I$, so faithfulness gives $g^p=1$.
For $h\in G$,
\[
(R_{g^h}-I)^2=R_h^{-1}(R_g-I)^2R_h=0,
\]
which proves the assertion about conjugates.
\end{proof}

\begin{lem}\label{lem:orbit}
Under the hypotheses of Theorem~\ref{thm:long}, put
$K=\langle a^{b^k}:k\in\mathbb Z\rangle$.
Then $K$ is elementary abelian of rank $r$ and is normalized by $b$.
Regard $K$ as an $\mathbb F_p$-vector space, with the identity
written as $0$. Let $T$ be conjugation by $b$ and put
$\Delta=T-I_K$. Then
\[
\Delta^i(a)=a_i\quad(0\leq i\leq r),\qquad T^p=I_K,
\]
and both
\[
a_0,\ldots,a_{r-1}\qquad\text{and}\qquad
w_0,\ldots,w_{r-1},\quad w_i=a^{b^i},
\]
are bases of $K$. Each $w_i$ is quadratic on $V$.
\end{lem}

\begin{proof}
The conjugates $a^{b^k}$ lie in the abelian subgroup $L$ and,
by Lemma~\ref{lem:quadratic}, have order $p$ and are quadratic.
Thus $K$ is elementary abelian and is normalized by $b$.
In additive notation,
$a_i=a_{i-1}^b-a_{i-1}=\Delta(a_{i-1})$, so
$a_i=\Delta^i(a)$. The hypotheses give
$\Delta^r(a)=0$ and $\Delta^{r-1}(a)\neq0$.
Since $r<p$, we have $T^p(a)=a+\Delta^p(a)=a$.
As $T^p$ commutes with $T$, it fixes every $T^k(a)$, and hence
$T^p=I_K$.

For $0\leq k<p$, the binomial theorem gives
\begin{equation}\label{eq:orbit-expansion}
T^k(a)=\sum_{j=0}^{r-1}\binom{k}{j}a_j.
\end{equation}
Since all $a_j$ belong to $K$, this shows that they span $K$.
If $\sum_{i=0}^{r-1}\lambda_i a_i=0$ is a nontrivial relation
and $j$ is the least index with $\lambda_j\neq0$, applying
$\Delta^{r-1-j}$ gives $\lambda_j a_{r-1}=0$, a contradiction.
Thus the $a_j$ form a basis. Finally,
$w_i=\sum_{j=0}^{i}\binom{i}{j}a_j$ for $0\leq i<r$,
so the change-of-basis matrix is triangular with diagonal entries
equal to $1$. Hence the $w_i$ form a quadratic basis of $K$.
\end{proof}

The next lemma linearizes the action of a group with a quadratic
basis by means of a truncated logarithm.

\begin{lem}\label{lem:log}
Let $K=\langle w_0,\ldots,w_{r-1}\rangle$
be an elementary abelian $p$-group of rank $r<p$,
regarded as an $\mathbb F_p$-vector space, where
$w_0,\ldots,w_{r-1}$ form a basis and each $w_i$
is quadratic on $V$.

Then there exists an $\mathbb F_p$-linear map
\[
\Phi:K\longrightarrow
\operatorname{End}_{\mathbb F_p}(V)
\]
such that
\[
\Phi(w)=\log R_w
\]
and
\[
(R_w-I)^2=0
\quad\Longleftrightarrow\quad
\Phi(w)^2=0
\]
for every $w\in K$.
\end{lem}

\begin{proof}
Put
\[
X_i=R_{w_i}-I,\qquad 0\leq i\leq r-1.
\]
Since the $w_i$ commute and are quadratic,
\[
X_iX_j=X_jX_i,\qquad X_i^2=0.
\]

Let
\[
\mathcal A=\mathbb F_p[I,X_0,\ldots,X_{r-1}],
\qquad
J=(X_0,\ldots,X_{r-1}).
\]
Then $\mathcal A$ is a commutative algebra.
Furthermore,
\[
J^{r+1}=0,
\]
because every monomial of degree $r+1$
in the $r$ commuting square-zero generators
contains some factor $X_i^2$.

Since $r<p$, we may define the truncated
logarithm and exponential by
\[
\log(I+U)
=\sum_{k=1}^{r}(-1)^{k+1}\frac{U^k}{k},
\]
and
\[
\exp(U)=\sum_{k=0}^{r}\frac{U^k}{k!}
\]
for $U\in J$.

These are the usual formal logarithm and exponential identities
truncated modulo terms of degree at least $r+1$. They are valid
because $J^{r+1}=0$ and all denominators that occur are invertible
in $\mathbb F_p$.

For any
\[
w=w_0^{\lambda_0}\cdots w_{r-1}^{\lambda_{r-1}},
\qquad \lambda_i\in\mathbb F_p,
\]
where each exponent is interpreted modulo $p$, we have
\[
R_w=\prod_{i=0}^{r-1}(I+\lambda_iX_i).
\]
Since $X_i^2=0$,
\[
I+\lambda_iX_i=\exp(\lambda_iX_i).
\]
Using the commutativity of the $X_i$, we obtain
\[
R_w
=\exp\left(\sum_{i=0}^{r-1}\lambda_iX_i\right).
\]
Consequently,
\[
\Phi(w)=\log R_w
=\sum_{i=0}^{r-1}\lambda_iX_i.
\]
This shows that $\Phi$ is $\mathbb F_p$-linear.

Finally, if $(R_w-I)^2=0$, then
\[
\Phi(w)=\log R_w=R_w-I,
\]
and hence $\Phi(w)^2=0$.

Conversely, if $\Phi(w)^2=0$, then
\[
R_w=\exp(\Phi(w))=I+\Phi(w),
\]
so
\[
(R_w-I)^2=\Phi(w)^2=0.
\]
The proof is complete.
\end{proof}

Applying this linearization to \eqref{eq:orbit-expansion} yields
an operator-valued polynomial to which interpolation applies.

\begin{lem}\label{lem:interpolation}
Under the hypotheses of Theorem~\ref{thm:long},
let $\Phi$ be the map in Lemma~\ref{lem:log},
and put
\[
A_j=\Phi(a_j),\qquad 0\leq j\leq r-1.
\]
Define
\[
F(t)=\sum_{j=0}^{r-1}\binom{t}{j}A_j,
\]
where
\[
\binom{t}{j}
=\frac{t(t-1)\cdots(t-j+1)}{j!}.
\]
Then
\[
F(t)^2\equiv0.
\]
In particular,
\[
A_{r-1}^2=0.
\]
\end{lem}

\begin{proof}
For every $t\in\mathbb F_p$, choosing an integer
representative of $t$, Lemma~\ref{lem:orbit} gives
\[
T^t(a)
=\sum_{j=0}^{r-1}\binom{t}{j}a_j.
\]
Applying the linear map $\Phi$, we obtain
\[
F(t)=\Phi(T^t(a)).
\]

Since
\[
T^t(a)=a^{b^t}
\]
is conjugate to $a$, it is quadratic on $V$.
By Lemma~\ref{lem:log},
\[
F(t)^2=0
\qquad\text{for every }t\in\mathbb F_p.
\]

Put $P(t)=F(t)^2$. Then $P(t)$ has degree at most
$2r-2\leq p-1$. Write
\[
P(t)=\sum_{k=0}^{p-1}C_k t^k,
\qquad C_k\in\mathcal A,
\]
where zero coefficients are inserted if necessary. Enumerate the
elements of $\mathbb F_p$ as $t_0,\ldots,t_{p-1}$. Since
$P(t_i)=0$ for every $i$, we have
\[
\begin{pmatrix}
1&t_0&\cdots&t_0^{p-1}\\
1&t_1&\cdots&t_1^{p-1}\\
\vdots&\vdots&&\vdots\\
1&t_{p-1}&\cdots&t_{p-1}^{p-1}
\end{pmatrix}
\begin{pmatrix}C_0\\C_1\\\vdots\\C_{p-1}\end{pmatrix}
=0.
\]
The determinant of this Vandermonde matrix is
\[
\prod_{0\leq i<j\leq p-1}(t_j-t_i),
\]
which is nonzero in $\mathbb F_p$. Hence the matrix is invertible
over $\mathbb F_p$. Its inverse acts on $\mathcal A^p$ as well,
so $C_0=\cdots=C_{p-1}=0$. Thus $P(t)\equiv0$.
This argument does not require $\mathcal A$ to be an integral domain.

Finally, the leading coefficient of $F(t)$ is
\[
\frac{A_{r-1}}{(r-1)!}.
\]
Therefore the coefficient of $t^{2r-2}$
in $F(t)^2$ is
\[
\frac{A_{r-1}^2}{((r-1)!)^2}.
\]
It follows that
\[
A_{r-1}^2=0.
\]
\end{proof}

\begin{proof}[Proof of Theorem~\ref{thm:long}]
By Lemma~\ref{lem:orbit}, $K=\langle a^{b^k}:k\in\mathbb Z\rangle$
has rank $r<p$ and a quadratic basis. Apply Lemma~\ref{lem:log}
to obtain $\Phi$. Lemma~\ref{lem:interpolation} gives
$\Phi(a_{r-1})^2=0$, so Lemma~\ref{lem:log} implies
$(R_{a_{r-1}}-I)^2=0$. Since $a_{r-1}\neq1$, it is quadratic.
\end{proof}

\begin{thm}\label{Z_i(G)}
Let $p$ be an odd prime, let $G$ be a finite $p$-group, and let
$V$ be a faithful $\mathbb F_pG$-module. If there is no quadratic
element in $\Omega_1(Z(G))$, then there is no quadratic element in
$Z_{(p+1)/2}(G)$.
\end{thm}

\begin{proof}
Put $m=(p+1)/2$. Suppose that $Z_m(G)$ contains a quadratic
element, and choose such an element $x$ for which $|C_G(x)|$ is
maximal. By Lemma~\ref{lem:quadratic}, $x$ has order $p$, so the
hypothesis implies that $x\notin Z(G)$. Set $M=C_G(x)$.
The normalizer condition for finite $p$-groups gives
\[
M<N_G(M).
\]
Choose $b\in N_G(M)\setminus M$ and put $L=Z(M)$. Then $L$ is
abelian, $x\in L$, and $b$ normalizes $L$, since $L$ is
characteristic in $M$.

Define $a_0=x$ and $a_i=[a_{i-1},b]$ for $i\geq1$. All these
elements belong to $L$. Also, $a_1\neq1$ because $b\notin C_G(x)$.
Since $x\in Z_m(G)$, we have $a_i\in Z_{m-i}(G)$ for
$0\leq i\leq m$, where $Z_0(G)=1$. In particular, $a_m=1$.
Thus there is an integer $r$ with $2\leq r\leq m$ such that
$a_{r-1}\neq1$ and $a_r=1$.
By Theorem~\ref{thm:long}, $a_{r-1}$ is quadratic on $V$.
Moreover, $a_{r-1}\in L\cap Z_m(G)$ and $[a_{r-1},b]=1$.
Consequently,
\[
M<\langle M,b\rangle\leq C_G(a_{r-1}),
\]
contradicting the maximal choice of $|C_G(x)|$.
\end{proof}

\begin{thm}\label{thm=C(G)x}
Let $p$ be an odd prime, let $G$ be a finite $p$-group, and let
$V$ be a faithful $\mathbb F_pG$-module. If there is no quadratic
element in $\Omega_1(Z(G))$, then $Z_{(p+1)/2}(G)$ centralizes
every quadratic element of $G$.
\end{thm}

\begin{proof}
Put $m=(p+1)/2$. Suppose that $x$ is quadratic on $V$ and that
$Z_m(G)\not\leq C_G(x)$. Set $M=C_G(x)$ and $H=MZ_m(G)$.
Since $Z_m(G)$ is normal in $G$, $H$ is a subgroup, and $M<H$.
The normalizer condition in $H$ gives $M<N_H(M)$.
Choose $y\in N_H(M)\setminus M$ and write $y=ub$ with
$u\in M$ and $b\in Z_m(G)$. Since $M\leq N_H(M)$, it follows
that $b=u^{-1}y\in N_H(M)\setminus M$.

Put $L=Z(M)$. Then $L$ is abelian, $x\in L$, and $b$ normalizes
$L$. Define $a_0=x$ and $a_i=[a_{i-1},b]$ for $i\geq1$.
We have $a_1\neq1$, while $b\in Z_m(G)$ implies
\[
a_i\in Z_{m-i}(G)\qquad(1\leq i\leq m).
\]
Hence $a_m=1$, and there is an integer $r$ with
$2\leq r\leq m$ such that $a_{r-1}\neq1$ and $a_r=1$.
Theorem~\ref{thm:long} shows that $a_{r-1}$ is quadratic.
Since $a_{r-1}\in Z_m(G)$, this contradicts
Theorem~\ref{Z_i(G)}.
\end{proof}

\section{Oliver's Conjecture and the Quadratic Conjecture}

We begin by recalling Oliver's $p$-group conjecture and its
module-theoretic formulation. We retain the left-normed
commutator convention from the Introduction. For subgroups
$A,B,C$, write $[A,B,C]=[[A,B],C]$, and, for $k\geq1$, let
$[A,B;k]=[A,\underbrace{B,\ldots,B}_{k\text{ copies}}]$.
For a subgroup $E\leq G$, put
\[
E^G=\langle E^g:g\in G\rangle,
\]
the normal closure of $E$ in $G$.

\begin{defi}[{\cite[Definition 2.1]{Lynd2013}}]
Let $3\leq k\leq p$. The $k$th Oliver subgroup of $S$ is the
largest subgroup $\mathfrak{X}_k(S)$ such that there exists a series
\[
1=Q_0<Q_1<\cdots<Q_{n-1}<Q_n=\mathfrak{X}_k(S)
\]
with each $Q_i$ normal in $S$ and with the property that for all
$1\leq i\leq n$,
\[
[\Omega_1(C_S(Q_{i-1})),Q_i;k-1]=1.
\]
\end{defi}

Taking $k=p$ in the definition, we get Oliver's original definition $\mathfrak{X}(S)$ in \cite{O}.
So $\mathfrak{X}_p(S)=\mathfrak{X}(S)$. Green et al. \cite{GHM2} also drew attention to
$\mathfrak{X}_3(S)$, which they denoted by $\mathcal{Y}(S)$.

Recall that the Thompson subgroup is the subgroup $J_e(S)$ of $S$
generated by the elementary abelian subgroups of maximum order.
In \cite[Conjecture 3.9]{O}, Oliver proposed the following conjecture.

\begin{conj}[Oliver]\label{conj:olive}
Let $S$ be a finite $p$-group with $p$ odd. Then
\[
J_e(S)\leq \mathfrak{X}(S).
\]
\end{conj}
In \cite[Conjecture 1.3]{GHM2}, Green et al. pose a strong form of Conjecture \ref{conj:olive} and its module-theoretic reformulation.
\begin{conj}[Strong Oliver Conjecture]\label{conj:strong}
Let $S$ be a finite $p$-group with $p$ odd. Then
\[
J_e(S)\leq \mathfrak{X}_3(S).
\]
\end{conj}
\begin{rem}\label{rem:strong-implies-oliver}
Since $\mathfrak{X}_3(S)\le \mathfrak{X}(S)$, if $S$ satisfies Conjecture~\ref{conj:strong}, then $S$ satisfies Conjecture~\ref{conj:olive}.
\end{rem}

To introduce the module-theoretic reformulation of Conjecture~\ref{conj:strong}, we recall the definitions of an $F$-module and an offender.
\begin{defi}[{\cite[26.5]{GLS}}]
Let $G$ be a finite group, and let $V$ be a faithful
finite-dimensional $\mathbb F_pG$-module. Suppose that a
nonidentity elementary abelian $p$-subgroup $E\leq G$ satisfies
\[
|E|\,|C_V(E)|\geq |V|,
\]
then $V$ is called an \emph{$F$-module} for $G$, and $E$
is called an \emph{offender}.
\end{defi}

Following Green, H\'ethelyi and Lilienthal, we say that a faithful
$\mathbb F_pG$-module $V$ satisfies property (PS) if its restriction
to every central subgroup of order $p$ has a nonzero projective
summand. Equivalently, every nonidentity element of
$\Omega_1(Z(G))$ acts on $V$ with minimal polynomial
$(X-1)^p$ \cite[Section~1]{GHL}.

The following reformulation of Oliver's conjecture will be used in
the proof of Corollary~\ref{thm:Oclass}.
\begin{thm}[Green--H\'ethelyi--Lilienthal, {\cite[Theorem~1.2]{GHL}}]
\label{thm:GHL-criterion}
Let $p$ be an odd prime and $G\neq1$ a finite $p$-group.
Oliver's $p$-group conjecture holds for every finite $p$-group $S$
with $S/\mathfrak X(S)\cong G$ if and only if $G$ has no faithful
$F$-modules satisfying (PS).
\end{thm}

For each subgroup $H\leq G$, put
\[
  j_H(V)=\frac{|H|\,|C_V(H)|}{|V|}.
\]
Thus a nontrivial elementary abelian $p$-subgroup $E$ is an
offender precisely when $j_E(V)\geq1$. An offender $A$ is called
a \emph{best offender} if
\[
j_A(V)\geq j_D(V)\qquad\text{for every }D\leq A.
\]
An offender $E$ is called a \emph{quadratic offender} if
$[V,E,E]=0$.

\begin{conj}[Quadratic Conjecture, {\cite[Conjecture~1.4]{GHM2}}]
\label{quadratic-conjecture}
Let $p$ be a prime and $G$ a finite $p$-group. If the faithful
$\mathbb F_pG$-module $V$ is an $F$-module, then there exists a quadratic element in $\Omega_1(Z(G))$.
\end{conj}

\begin{rem}\label{rem:strong-module-equivalence}
By Green, H\'ethelyi and Mazza \cite[Theorem~1.5]{GHM2},
for an odd prime $p$ and a finite $p$-group $G$, the Strong
Oliver Conjecture holds
for every finite $p$-group $S$ with
$S/\mathfrak X_3(S)\cong G$ if and only if the Quadratic
Conjecture holds for every faithful $\mathbb F_pG$-module.
\end{rem}

We use the following formulation of Timmesfeld's replacement
theorem, including the explicit construction given in
\cite[Lemma~4.1 and its proof]{GHM2010}.

\begin{lem}[Timmesfeld's replacement theorem, {\cite[Lemma~4.1 and its proof]{GHM2010}}]
\label{lem:replacement}
Let $G$ be a finite $p$-group and let $V$ be a finite-dimensional
faithful $\mathbb{F}_pG$-module. Suppose that $A$ is a nontrivial
elementary abelian subgroup of $G$ satisfying
\[
    j_A(V)\geq j_D(V)
    \quad\text{for all }D\leq A.
\]
Set
\[
    B=C_A([V,A]).
\]
Then $B$ is a quadratic offender contained in $A$.
In particular,
\[
    j_B(V)=j_A(V)
    \quad\text{and}\quad
    [V,B,B]=0.
\]
\end{lem}

Immediately we have the following theorem.

\begin{thm}\label{thm:smallclass}
Let $p$ be an odd prime and $G$ a finite $p$-group. If the
nilpotency class of $G$ is at most $\frac{p+1}{2}$, then $G$
satisfies the Quadratic Conjecture.
\end{thm}
\begin{proof}
Put $m=(p+1)/2$, so that $Z_m(G)=G$. Let $V$ be a faithful
$F$-module for $G$. Among all offenders, choose $A$ for which
$j_A(V)$ is maximal. If $D\leq A$ and $D\neq1$, then either
$D$ is an offender, in which case $j_D(V)\leq j_A(V)$, or
$j_D(V)<1\leq j_A(V)$. Also $j_1(V)=1\leq j_A(V)$.
Thus $A$ is a best offender.
By Lemma~\ref{lem:replacement}, $A$ contains a nontrivial
quadratic offender $B$. Therefore $G=Z_m(G)$ contains a
quadratic element. The contrapositive of Theorem~\ref{Z_i(G)}
now gives a quadratic element in $\Omega_1(Z(G))$.
\end{proof}

We use the weak-closure convention from the Introduction. The next
lemma converts a commutator condition into normalization.

\begin{lem}\label{lem:weakly closed}
Suppose that $E$ is a weakly closed abelian subgroup of $G$.
If $H\leq G$ satisfies $[H,E,E]=1$, then $H\leq N_G(E)$.
\end{lem}

\begin{proof}
Let $h\in H$ and $e\in E$. The hypothesis gives
$[h,e]\in C_G(E)$, and hence $[e,h]=[h,e]^{-1}\in C_G(E)$.
Since $E$ is abelian, $e^h=e[e,h]\in C_G(E)$.
Thus $[E,E^h]=1$. Weak closure forces $E^h=E$, so
$h\in N_G(E)$. This holds for every $h\in H$.
\end{proof}

The following theorem provides a weakly closed quadratic offender whenever $V$ is an $F$-module.

\begin{thm}[{\cite[Proposition~4.5]{GHM2}}]\label{weakly}
Suppose that a faithful $\mathbb{F}_p[G]$-module $V$ is an F-module. Set
$$j_0 =\max\{ j_E(V)\mid E \ \mbox{an offender}\}.$$
Then there is a weakly closed quadratic offender $E$ with $j_E(V) = j_0$.
Moreover if $D\le G$ is any offender with $j_D(V)=j_0$, then there is such an $E$ which is a subgroup of the
normal closure of $D$.
\end{thm}

\begin{thm}\label{thm=p+3/2}
Let $p$ be an odd prime, let $G$ be a finite $p$-group, and let
$V$ be a faithful $\mathbb F_pG$-module. Suppose that $E$ is a
nontrivial weakly closed elementary abelian subgroup with
$[V,E,E]=0$. If there is no quadratic element
in $\Omega_1(Z(G))$, then $[Z_{(p+3)/2}(G),E]=1$.
\end{thm}

\begin{proof}
Put $m=(p+1)/2$ and $H=Z_{m+1}(G)$. Since $[H,E]\leq Z_m(G)$,
Theorem~\ref{thm=C(G)x} implies $[H,E,E]=1$: every nonidentity
element of $E$ is quadratic. Lemma~\ref{lem:weakly closed}
therefore gives $H\leq N_G(E)$. It follows that
\[
[H,E]\leq E\cap Z_m(G).
\]
Every nonidentity element of this intersection would be quadratic,
contrary to Theorem~\ref{Z_i(G)}. Hence $[H,E]=1$.
\end{proof}

\section{Relative lower central series}

Let $G$ be a finite $p$-group and $N\unlhd G$. The
\emph{relative lower central series} of $N$ in $G$ is defined by
\[
\Gamma_1(N)=N,\qquad
\Gamma_{i+1}(N)=[\Gamma_i(N),G]\quad(i\geq1).
\]
See \cite[Definition~1.1]{Hilton1976}. The ambient group $G$ is
fixed throughout, so it is suppressed from the notation.
Each $\Gamma_i(N)$ is normal in $G$, and the series is descending.
Recall also the lower central series of $G$,
\[
\gamma_1(G)=G,\qquad
\gamma_{i+1}(G)=[\gamma_i(G),G]\quad(i\geq1).
\]
Induction gives $\Gamma_i(N)\leq\gamma_i(G)$ for every $i$.
Consequently, if $G$ has nilpotency class $c$, then
$\Gamma_{c+1}(N)=1$.

\begin{lem}\label{lem:normal closure}
Let $E$ be a subgroup of a group $G$, put $N=E^G$, and let
$H\unlhd G$. If $B=[H,E]$, then
\[
[H,N]=B^G=\langle B^g:g\in G\rangle.
\]
\end{lem}

\begin{proof}
Put $D=B^G$. Since $H$ and $N$ are normal in $G$, $[H,N]$
is normal in $G$. It contains $B$, so $D\leq[H,N]$.
For the reverse inclusion, work in $G/D$ and denote images by
bars. By construction, $\overline H$ centralizes $\overline E$.
As $\overline H$ is normal in $\overline G$, it also centralizes
every conjugate of $\overline E$. Hence it centralizes
$\overline N=\langle\overline E^{\,\overline g}:
\overline g\in\overline G\rangle$. Thus $[H,N]\leq D$.
\end{proof}

The next lemma justifies the commutator interchange needed below.
For convenience, set $[A,G;0]=A$.

\begin{lem}\label{lem:commutator-shifting}
Let $H\leq N$ be normal subgroups of a group $G$, and suppose
that $[H,G]\leq Z(N)$. Then, for every integer $k\geq0$,
\[
[[H,N],G;k]=[H,\Gamma_{k+1}(N)].
\]
\end{lem}

\begin{proof}
We first note that if $A,B\unlhd G$ and $[[A,G],B]=1$, then
\begin{equation}\label{eq:commutator-interchange}
[[A,B],G]=[A,[B,G]].
\end{equation}
Indeed, put $P=[[A,B],G]$ and $Q=[A,[B,G]]$; both are normal
in $G$. In $G/Q$, the subgroups $[[B,G],A]$ and $[[G,A],B]$
are trivial, so the Three Subgroups Lemma gives $P\leq Q$.
Applying the same lemma in $G/P$ gives $Q\leq P$.
This proves \eqref{eq:commutator-interchange}.

We now induct on $k$. The case $k=0$ is immediate. If the
assertion holds for $k$, then $\Gamma_{k+1}(N)\leq N$ and
$[H,G]\leq Z(N)$ imply
\[
[[H,G],\Gamma_{k+1}(N)]=1.
\]
Applying \eqref{eq:commutator-interchange}, we obtain
\begin{align*}
[[H,N],G;k+1]
&=[[H,\Gamma_{k+1}(N)],G]\\
&=[H,[\Gamma_{k+1}(N),G]]\\
&=[H,\Gamma_{k+2}(N)],
\end{align*}
as required.
\end{proof}

\begin{lem}\label{lem:weakly close-series}
Let $G$ be a finite $p$-group, and let $1\neq E\leq G$ be a
weakly closed abelian subgroup. Suppose that $N=E^G$ is
nonabelian, and let $t$ be the least integer such that
$\Gamma_t(N)\leq Z(N)$. Then $t\geq2$,
$\Gamma_{t-1}(N)\leq N_G(E)$, and, setting
$B=[\Gamma_{t-1}(N),E]$, we have
\[
1\neq B\leq E,\qquad
1\neq[\Gamma_{t-1}(N),N]=B^G.
\]
\end{lem}

\begin{proof}
Since $N$ is nonabelian, $t\geq2$. Put $H=\Gamma_{t-1}(N)$.
By the definition of $t$,
\[
[H,E]\leq[H,G]=\Gamma_t(N)\leq Z(N).
\]
As $E\leq N$, this gives $[H,E,E]=1$.
Lemma~\ref{lem:weakly closed} now implies $H\leq N_G(E)$, so
$B=[H,E]\leq E$.
Since $H\unlhd G$, Lemma~\ref{lem:normal closure} gives
$[H,N]=B^G$. The minimality of $t$ implies that
$H\not\leq Z(N)$. Since $H\leq N$, we have $[H,N]\neq1$,
and hence $B\neq1$ as well.
\end{proof}

\begin{proof}[Proof of Theorem~\ref{thm:weakly-class}]
Put $m=(p+1)/2$ and $c=\operatorname{cl}(G)$, so that
$c\leq p+2=2m+1$. Let $V$ and $E$ be as in the statement,
and suppose, for a contradiction, that $\Omega_1(Z(G))$
contains no quadratic element.

Put $N=E^G$. We first show that $N$
is nonabelian. If $N$ were abelian, then $[E,E^g]=1$ for every
$g\in G$, so weak closure would imply $E^g=E$ for every $g$.
Thus $E$ would be a nontrivial normal subgroup of the finite
$p$-group $G$, and hence $E\cap Z(G)\neq1$. Since $E$ is
elementary abelian and quadratic, every nonidentity element of
this intersection would be a quadratic element of $\Omega_1(Z(G))$,
a contradiction.

Let $t$ be the least integer such that $\Gamma_t(N)\leq Z(N)$.
Since $N$ is nonabelian and $\Gamma_{c+1}(N)=1$, we have
$2\leq t\leq c+1$. Set
$H=\Gamma_{t-1}(N)$. We claim that $t\leq m+1$.
Since $\Gamma_{t-1}(N)\leq\gamma_{t-1}(G)$ and
$\gamma_{t-1}(G)\leq Z_{c-t+2}(G)$, we have
\[
H\leq Z_{c-t+2}(G).
\]
Here the latter standard containment follows directly from
$[\gamma_{t-1}(G),G;c-t+2]=\gamma_{c+1}(G)=1$.
If $t\geq m+2$, then
\[
c-t+2\leq(2m+1)-(m+2)+2=m+1,
\]
so $H\leq Z_{m+1}(G)$. By Theorem~\ref{thm=p+3/2},
$H$ centralizes $E$. Since $H\unlhd G$,
Lemma~\ref{lem:normal closure} implies $[H,N]=1$.
As $H\leq N$, this says that $H\leq Z(N)$, contrary to the
minimality of $t$. Therefore $t\leq m+1$.

By Lemma~\ref{lem:weakly close-series},
\[
D:=[H,N]=B^G\neq1,\qquad 1\neq B:=[H,E]\leq E.
\]
Also $[H,G]=\Gamma_t(N)\leq Z(N)$, so
Lemma~\ref{lem:commutator-shifting} gives
\[
[D,G;t-1]=[H,\Gamma_t(N)]=1,
\]
where the last equality follows from $H\leq N$ and
$\Gamma_t(N)\leq Z(N)$. Consequently,
\[
B\leq D\leq Z_{t-1}(G)\leq Z_m(G).
\]
Choose $1\neq b\in B$. Since $B\leq E$ and $E$ is quadratic,
$b$ is quadratic on $V$. This contradicts
Theorem~\ref{Z_i(G)}. The contradiction proves the theorem.
\end{proof}

\begin{proof}[Proof of Corollary~\ref{thm:class}]
Let $V$ be a faithful $F$-module for $G$.
By Theorem~\ref{weakly}, $G$ contains a nontrivial weakly closed
quadratic offender $E$. In particular, $E$ is elementary abelian
and $[V,E,E]=0$. Theorem~\ref{thm:weakly-class} therefore gives
a quadratic element in $\Omega_1(Z(G))$.
\end{proof}

\begin{proof}[Proof of Corollary~\ref{thm:Oclass}]
Let $G=S/\mathfrak X(S)$. If $G=1$, then
$\mathfrak X(S)=S$ and Oliver's conjecture holds trivially.
Suppose that $G\neq1$. By hypothesis,
$\operatorname{cl}(G)\leq p+2$.
By Corollary~\ref{thm:class}, every faithful $F$-module
$V$ for $G$ contains a nonidentity quadratic element
$z\in\Omega_1(Z(G))$. Since $V$ is faithful, the minimal
polynomial of $z$ is exactly $(X-1)^2$.
On the other hand, if $V$ satisfied (PS), then $z$ would act
with minimal polynomial $(X-1)^p$. Since $p$ is odd, these
conditions are incompatible. Hence $G$ has no faithful
$F$-modules satisfying (PS). Theorem~\ref{thm:GHL-criterion}
therefore implies that $S$ satisfies Oliver's $p$-group
conjecture.
\end{proof}

\section*{Financial Support}
This work was supported by NSFC (Nos.~12371022 and 12271318).

\section*{Acknowledgments}
During the preparation of this manuscript, ChatGPT-6 Extra High(OpenAI) was used for assistance with literature
searches, discussion of proof strategies, language editing, and LaTeX
preparation. ChatGPT's most significant contribution to the mathematical development of this work was giving and proving Lemma \ref{lem:log}.
The authors are responsible for verifying all AI-assisted
suggestions and take full responsibility for the accuracy, originality,
and integrity of the final manuscript, including its mathematical
statements, proofs, and references.

\end{document}